\documentclass[a4paper,12pt]{article}

\usepackage[utf8]{inputenc}
\usepackage[mathscr]{eucal}
\usepackage{amsmath}
\usepackage{amssymb}
\usepackage{MnSymbol}
\usepackage{amsfonts}
\usepackage{graphicx}
\usepackage{multicol}
\usepackage{enumerate}

\newcommand{\zbp}{\emptyset}

\newcommand{\Q}{\mathbb{Q}}
\newcommand{\vare}{\varepsilon}
\newcommand{\varp}{\varphi}
\newcommand{\Fin}{\text{Fin}}
\newcommand{\oF}{{\mathcal{F}}}
\newcommand{\cont}{{\mathfrak{c}}}
\newcommand{\bont}{{\mathfrak{b}}}
\newcommand{\dont}{{\mathfrak{d}}}
\newcommand{\NN}{\boldsymbol{\text{N}}}
\newcommand{\non}{\text{non}}

\mathchardef\mhyphen="2D

\newtheorem{theo}{Theorem}
\newtheorem{lemma}[theo]{Lemma}

\newtheorem{cor}[theo]{Corollary}

\title{Generalized Egorov's statement for ideals}
\author{Michał Korch}
\date{}

\begin{document}

\maketitle

\begin{abstract}
We consider the generalized Egorov's statement (Egorov's Theorem without the assumption on measurability of the functions, see \cite{tw:nget}) in the case of an ideal convergence and a number of different types of ideal convergence notion. We prove that in those cases the generalized Egorov's statement is independent from ZFC.
\end{abstract}

\section{Introduction}

In this paper we consider various versions of the classic Egorov's Theorem. Let us recall (see e.g. \cite{jo:mc}) that the classic Egorov's Theorem states  that given a sequence of measurable functions (we restrict our attention to the real functions $[0,1]\to[0,1]$) which is pointwise convergent on $[0,1]$ and $\varepsilon>0$, one can find a measurable set $A\subseteq [0,1]$ with $m(A)\geq 1-\vare$ such that the sequence converges uniformly on $A$ ($m$ denotes the Lebesgue measure).

It is interesting whether we can drop the assumption on measurability of the functions in the above theorem. A statement which says that given any sequence of functions $[0,1]\to[0,1]$ which is pointwise convergent and $\vare>0$, there exists a set $A\subseteq [0,1]$ with $m^*(A)\geq 1-\vare$ ($m^*$ denotes the outer measure) such that the sequence converges uniformly on $A$, is called the generalized Egorov's statement. T. Weiss in his unpublished manuscript (see \cite{tw:nget}) proved that it is independent from ZFC, and this fact was used in \cite{fdasostw:ssa}. Then R.~Pinciroli studied the method of T.~Weiss more systematically (see \cite{rp:igsetzw}). For example, he related it to cardinal coefficients: $\non(\NN)$ (the lowest possible cardinality of a non-null set), $\bont$ (the lowest possible cardinality of a family of sequences of natural numbers unbounded in the sense of the order $\leq^*$ of eventual domination) and $\dont$ (the lowest possible cardinality of a family of sequences of natural numbers such that every possible sequence is dominated in the sense of $\leq^*$ by a sequence in the family). In particular, he proved that $\non(\NN)<\mathfrak{b}$ implies that the generalized Egorov's statement holds, but if, for example, $\non(\NN)=\dont=\cont$, then it fails.

We can also define a notion of convergence of a sequence of functions with respect to a given ideal $I$ on $\omega$. There are different types of convergence with respect to $I$, and pointwise and uniform convergence are the most common. Given two notion of convergence with respect to an ideal, we can ask whether the classic Egorov's Theorem (with the measurability assumption) holds for those two notion of convergence in the sense of whether the weaker convergence implies the stronger convergence on a subset of arbitrarily large measure. The answer may often be negative as in the case of uniform and pointwise convergence for many analytic P-ideals (see \cite[Theorem 3.4]{nm:ivetap}). But one can also define other types of convergence, e.g. equi-ideal convergence. And, for example, in the case of analytic P-ideal so called weak Egorov's Theorem for ideals (between equi-ideal and pointwise ideal convergence) was proved by N.~Mro\.{z}ek (see \cite[Theorem 3.1]{nm:ivetap}). 

Therefore, we ask whether in the case of an ideal and two notion of convergence for which the Egorov's theorem with measurability assumption holds, we can drop this assumption. This paper deals with this question in relation to different types of ideal convergence notion.

\section{Using Pinciroli's method}

We start by a generalization of the method presented by R.~Pinciroli (see \cite{rp:igsetzw}, and also \cite{mr:get}). The core of this method can be generalized to the following theorem.

\begin{theo}\label{thpos1}
Assume that $\non(\NN)<\bont$. Let $\Phi \in (\omega^{\omega})^{[0,1]}$. Then for any $\vare>0$, there exists $A\subseteq [0,1]$ such that $m^{*}(A)\geq 1-\vare$ and $\Phi$ is bounded on $A$.
\end{theo}

Proof: We follow the arguments of Pinciroli (see \cite{rp:igsetzw}). 

Assume that $\non(\NN)<\bont$. Notice that this statement holds for example in a model obtained by $\aleph_2$-iteration with countable support of Laver forcing (see e.g. \cite{tbhj:stsrl}). Also it can be easily proven, that under this assumption there exists a set $Y\subseteq [0,1]$ of cardinality less that $\bont$ such that $m^{*}(Y)=1$. Indeed, if $N\subseteq [0,1]$ is a set of positive outer measure with $|N|<\bont$, then let $Y=\{x+y\colon x\in N, y\in \Q\}$, where $+$ denotes addition modulo $1$. Then $Y$ has outer measure $1$ under the Zero-One Law.

Therefore, every function $\varp\colon [0,1]\to\omega^\omega$ maps $Y$ onto a $K_{\sigma}$-set, where $K_{\sigma}$ denotes the $\sigma$-ideal of subsets of $\omega^\omega$ generated by the compact (equivalently bounded) sets. We get that $\Phi[Y]\in K_{\sigma}$. Assume that $\Phi[Y]\subseteq \bigcup_{n\in\omega} B_n$ with each $B_{n}$ bounded. Let $A_{n}=\Phi^{-1}[\bigcup_{i=0}^n B_{i}]$. Therefore,  $\Phi[A_n]$ is bounded, and for any $\vare>0$, there exists $n\in\omega$ such that $m^{*}(A_{n})\geq 1-\vare$. \hfill $\square$

In the products of the form $\omega^S$ and $(\omega^S)^T$ we consider
the partial orderings, denoted by the same symbol~$\leq$, given by
$x\leq y$, if $x(s)\leq y(s)$ for $x,y\in\omega^S$, $s\in S$,
and $\phi\le\psi$, if $\phi(t)\le\phi(t)$ for $\phi,\psi\in(\omega^S)^T$,
where $\phi(t),\psi(t)\in\omega^S$.
We say that a~function $o:X\to P$ from a~set~$X$ into a~partially ordered
set~$P$ is cofinal if for every $p\in P$ there exists $x\in X$ such that
$p\le o(x)$.

For a~sequence of functions $f_n:[0,1]\to[0,1]$ and subsets $A\subseteq[0,1]$
we consider notion of convergence $f_n\looparrowright f$ on~$A$.
We assume that if $B\subseteq A$ and $f_n\looparrowright f$ on~$A$,
then $f_n\looparrowright f$ on~$B$.
We write $f_n\looparrowright f$ provided that $f_n\looparrowright f$ on~$[0,1]$. Let $\oF\subseteq \{\langle f_n\rangle_{n\in\omega}:\forall_{n\in\omega} f_n:[0,1]\to[0,1]\}$ be an arbitrary family of sequences of functions. 

We consider two hypotheses between $\oF$ and $\looparrowright$:
\begin{description}
\item[$(H^\Rightarrow(\oF,\looparrowright))$]
There exists $o:\oF\to(\omega^\omega)^{[0,1]}$ such that for
every $F\in\oF$ and every $A\subseteq[0,1]$, if $o(F)[A]$
is bounded in $(\omega^\omega,\leq)$, then $F\looparrowright 0$ on~$A$.
\item[$(H^\Leftarrow(\oF,\looparrowright))$]
There exists cofinal $o:\oF\to(\omega^\omega)^{[0,1]}$
such that for every $F\in\oF$ and every $A\subseteq[0,1]$,
if $F\looparrowright0$ on~$A$, then $o(F)[A]$ is bounded in
$(\omega^\omega,\leq)$.
\end{description}

\begin{theo}\label{thpos}
Assume that $\non(\NN)<\bont$, and $H^\Rightarrow(\oF, \looparrowright)$. Then for any $\left<f_{n}\right>_{n\in\omega}\in\oF$ and any $\vare>0$, there exists $A\subseteq [0,1]$ such that $m^{*}(A)\geq 1-\vare$ and $f_{n}\looparrowright 0$ on $A$.
\end{theo}

Proof: Apply Theorem \ref{thpos1} for $o(\left<f_{n}\right>_{n\in\omega})$ given by $H^\Rightarrow(\oF, \looparrowright)$. \hfill $\square$

Recall that $Z\subseteq \omega^{\omega}$ is a $\cont$-Lusin set if it is of cardinality $\cont$, and if $A\subseteq Z$ is meagre, then $|A|<\cont$. The existence of such a set is independent from ZFC. Notice also that there exists a model of ZFC in which $\non(\NN)=\cont$, and there exists $\cont$-Lusin set. To get this model it suffices to iterate $\aleph_2$-times Cohen forcing with finite supports over a model of GCH (see \cite[Model~7.5.8 and Lemma~8.2.6]{tbhj:stsrl}).

\begin{theo}\label{thneg}
Assume that $\non (\NN)=\cont$, and that there exists a $\cont$-Lusin set. If $H^\Leftarrow(\oF,\looparrowright)$ holds, then there exist $\left<f_{n}\right>_{n\in\omega}\in\oF$ and $\vare>0$ such that for all $A\subseteq [0,1]$ with $m^{*}(A)\geq 1-\vare$, $f_{n}\not\looparrowright 0$ on $A$.
\end{theo}

Proof: Again, we follow the arguments of Pinciroli (see \cite{rp:igsetzw}). Let $Z\subseteq \omega^{\omega}$ be a $\cont$-Lusin set. Since every compact set is meagre in $\omega^\omega$, every $K_\sigma$ set is also meagre. Therefore, if $A\subseteq Z$ is a $K_\sigma$ set, then $|A|<\cont$. Let $o:\oF\to(\omega^\omega)^{[0,1]}$ be a cofinal function given by $H^\Leftarrow(\oF,\looparrowright)$. Let $\varp$ be a bijection between $[0,1]$ and $Z$. Finally, let $\left<f_{n}\right>_{n\in\omega}=F\in\oF$ be such that $o(F)\geq \varp$. 

To get a contradiction, assume that for every $i\in\omega$, there exists $A_i\subseteq [0,1]$ such that $m^*(A_i)\geq 1-1/2^i$ and $f_{n}\looparrowright 0$ on $A_i$. Let $A=\bigcup_{i\in\omega}A_i$.
For any $i\in\omega$, $o(F)[A_i]$ is bounded because $f_n\looparrowright 0$
on~$A_i$, and so $\varphi[A_i]$ is bounded since $o(F)\geq\varphi$.
Therefore, $\varphi[A]\in K_\sigma$ and $|A|=|\varphi[A]|<\cont$ because
$\varphi[A]\subseteq Z$.
This is a~contradiction because $m^*(A)=1$ and $\non(\NN)=\cont$. \hfill $\square$

The following theorem was proved by R.~Pinciroli in \cite{rp:igsetzw}.

\begin{cor}\label{oryg}
Assume that $\non(\NN)<\bont$. Then  for any $\left<f_{n}\right>_{n\in\omega}$ such that $f_n\colon[0,1]\to[0,1]$ for $n\in\omega$, and $f_n\to 0$, and any $\vare>0$, there exists $A\subseteq [0,1]$ such that $m^{*}(A)\geq 1-\vare$ and $f_{n}\rightrightarrows 0$ on $A$.\\
On the other hand, assume that $\non (\NN)=\cont$, and that there exists a $\cont$-Lusin set. Then there exist $\left<f_{n}\right>_{n\in\omega}$ such that $f_n\colon [0,1]\to [0,1]$ for $n\in\omega$, and $f_n\to 0$, and $\vare>0$ such that for all $A\subseteq [0,1]$ with $m^{*}(A)\geq 1-\vare$, $f_{n}\not\rightrightarrows 0$ on $A$.
\end{cor}

Proof: Let $\left<f_{n}\right>_{n\in\omega}$ be such that $f_{n}\to 0$. Set $\vare_{n}=1/2^n$, $n\in\omega$. Consider $\oF=\{\langle f_n\rangle_{n\in\omega}:\forall_{n\in\omega} f_n:[0,1]\to[0,1]\land f_n\to 0\}$ and $\looparrowright=\rightrightarrows$. Define $o:\oF\to(\omega^\omega)^{[0,1]}$ in the following way. Let \[o F(x)(n)=\min\{m\in\omega\colon \forall_{l\geq m} f_{l}(x)\leq\vare_n\}.\] We  get exactly the reasoning and the results of R.~Pinciroli (see \cite{rp:igsetzw}). He proves that the above function $o$ proves that both $H^\Leftarrow(\oF_\rightarrow,\rightrightarrows)$ and $H^\Rightarrow(\oF_\rightarrow,\rightrightarrows)$ hold, and then proves Theorems~\ref{thpos} and \ref{thneg} in this particular case. \hfill $\square$

In next sections we apply the method used in the proof of Corollary~\ref{oryg}. Assume that we are given two notion of convergence of sequences of functions
$f_n\leadsto f$ and $f_n\looparrowright f$ such that
$f_n\looparrowright f$ implies $f_n\leadsto f$.
We take
\[
\oF_{\leadsto}=\{\langle f_n\rangle_{n\in\omega}\colon
\forall_{n\in\omega}\ f_n:[0,1]\to[0,1]\land f_n\leadsto 0\}
\]
and we apply Theorem~\ref{thpos} and Theorem~\ref{thneg} with a~suitable function
$o:\oF_{\leadsto}\to(\omega^\omega)^{[0,1]}$ to get a~conclusion on the stronger
convergence $f_n\looparrowright 0$ of sequences from~$\oF_{\leadsto}$.

\section{Pointwise and equi-ideal convergence (for analytic $P$-ideals)}

Let $I$ be an analytic $P$-ideal and $f_n\colon[0,1]\to [0,1]$, $n\in\omega$. 
By the well-known result of Solecki $I=Exh(\phi)$ (\cite{ss:aita}), where $\phi$ is a lower semicontinuous submeasure (a function $\phi\colon 2^{\omega}\to [0,\infty]$ satisfying the following conditions: $\phi(\zbp)=0, \phi(A)\leq\phi (A\cup B)\leq \phi(A)+\phi(B)$ and $\phi(A)=\lim_{n\to\omega}\phi(A\cap n)$, for any $A,B\subseteq \omega$) and $Exh(\phi)=\{A\subseteq \omega\colon \lim_{n\to\infty}\phi(A\setminus n)=0\}$ (see also \cite{nm:ivetap}).

Fix a lower continuous submeasure $\phi$ such that $I=Exh(\phi)$. Recall that we have the following notion of convergence (see \cite{nm:ivetap}) on a set $A\subseteq [0,1]$:

\begin{description}
\item[pointwise ideal, $f_{n}\to_{I} 0$]
if and only if
  \[\forall_{\vare>0}\forall_{x\in A}\exists_{k\in\omega}\phi(\{n\in\omega\colon f_{n}(x)\geq\vare\}\setminus k)<\vare,\]
\item [equi-ideal, $f_{n}\twoheadrightarrow_{I} 0$]
if and only if \[\forall_{\vare>0}\exists_{k\in\omega}\forall_{x\in A}\phi(\{n\in\omega\colon f_{n}(x)\geq\vare\}\setminus k)<\vare,\]
\item[uniform ideal, $f_{n}\rightrightarrows_{I} 0$] if and only if \[\forall_{\vare>0}\exists_{k\in\omega}\phi(\{n\in\omega\colon\sup_{x\in A} f_{n}(x)\geq\vare\}\setminus k)<\vare.\]
\end{description}

It was proved in \cite{nm:ivetap} that these notion of convergence are independent from
the submeasure representation of~$I$. Moreover, the pointwise ideal and uniform ideal convergences can be expressed without the notion of a~submeasure and they coincide with the notion of
well-known ideal convergences defined for any ideal~$I$ on~$\omega$
(see the next section and also~\cite{nm:zicf}).

Obviously, $f_{n}\rightrightarrows_{I} 0\Rightarrow f_{n}\twoheadrightarrow_{I} 0\Rightarrow f_{n}\to_{I} 0$.

It was also proved in \cite{nm:ivetap} that the ideal version of Egorov's Theorem holds (in the case of analytic $P$-ideals) between equi-ideal and pointwise ideal convergence, i.e. if $\left<f_n\right>_{n\in\omega}$ is a sequence of measurable functions with $f_{n}\to_{I} 0$ on $[0,1]$ and $\vare>0$, then there exists $A\subseteq [0,1]$ such that $m(A)\geq 1-\vare$ and $f_{n}\twoheadrightarrow_{I} 0$ on $A$. Moreover, it was proved that the ideal version of Egorov's Theorem (in the case of analytic $P$-ideals) does not hold between uniform ideal and pointwise ideal convergence except for the trivial and ``pathological'' cases (see also \cite{nm:zicf}).

Notice that since $I$ is a proper ideal, $\lim_{i\to\infty} \phi(\omega\setminus i)>0$. If $\lim_{i\to\infty} \phi(\omega\setminus i)<\infty$, let \[\vare_n=\frac{\lim_{i\to\infty} \phi(\omega\setminus i)}{2^{n+1}}\] for $n\in\omega$. Otherwise set $\vare_n=1/2^{n+1}$. To use the method described in the previous section, we state the following definition. For a sequence of functions $F=\left<f_{n}\right>_{n\in\omega}, f_{n}\colon [0,1]\to [0,1]$ such that $f_{n}\to_{I} 0$, let $o_{\phi} F\in (\omega^{\omega})^{[0,1]}$, and
\[(o_{\phi} F)(x)(n)=\min\{k\in\omega\colon \phi(\{m\in\omega\colon f_{m}(x)\geq\vare_{n}\}\setminus k)<\vare_{n}\}.\]

The function $o_{\phi}\colon\oF_{\to_I}\to (\omega^\omega)^{[0,1]}$ is well defined, because for each $n\in\omega$, $\{k\in\omega\colon \phi(\{m\in\omega\colon f_{m}(x)\geq\vare_{n}\}\setminus k)<\vare_{n}\}$ is not empty since $f_{n}\to_{I} 0$.

\begin{lemma}\label{eqlem}
Let $F=\left<f_{n}\right>_{n\in\omega}$ be a sequence of functions such that $f_n\colon [0,1]\to [0,1]$. Then $f_{n}\twoheadrightarrow_{I} 0$ on $A\subseteq [0,1]$ if and only if $ (o_\phi(\left<f_{n}\right>_{n\in\omega}))[A]$ is bounded in $\omega^\omega$. In particular, $H^\Rightarrow(\oF_{\to_I},\twoheadrightarrow_I)$ holds.
\end{lemma}
Proof: By definition, $f_{n}\twoheadrightarrow_{I} 0$ on $A$ if and only if for any $n\in\omega$, there exists $k\in\omega$ such that for all $x\in A$, $\phi(\{m\in\omega\colon f_{m}(x)\geq\vare_{n}\}\setminus k)<\vare_{n}$. This is true if and only if there exists a sequence $\left<k_{n}\right>_{n\in\omega}$ of natural numbers such that for any $n\in\omega$ and $x\in A$, $\phi(\{m\in\omega\colon f_{m}(x)\geq\vare_{n}\}\setminus k_{n})<\vare_{n}$, which holds if and only if for all $x\in A$, $(o_{\phi} F) (x)(n)\leq k_{n}$. \hfill $\square$

\begin{cor}
Assume that $\non(\NN)<\bont$. Let $I$ be any analytic $P$-ideal, $\vare>0$, and let $F=\left<f_{n}\right>_{n\in\omega}$, $f_{n}\colon [0,1]\to [0,1]$  for $n\in\omega$, be such that $f_{n}\to_{I} 0$. Then there exists $A\subseteq [0,1]$ with $m^{*}(A)\geq 1-\vare$ such that $f_{n}\twoheadrightarrow_{I} 0$ on $A$ (the ideal version of the generalized Egorov's statement between equi-ideal and pointwise ideal convergence for analytic $P$-ideals is consistent with ZFC).
\end{cor}
Proof: Apply Theorem \ref{thpos} and Lemma \ref{eqlem}. \hfill $\square$

\begin{lemma}\label{organicz}
For any $\varp\colon [0,1]\to \omega^{\omega}$, there exists $F=\left<f_{n}\right>_{n\in\omega}$, $f_{n}\colon [0,1]\to [0,1]$ for $n\in\omega$ with $f_{n}\to_{I} 0$  such that $o_{\phi} F\geq\varp$. In particular, $H^\Leftarrow(\oF_{\to_I},\twoheadrightarrow_I)$ holds.
\end{lemma}
Proof: Fix $x\in [0,1]$. Notice that $\phi(\omega\setminus n)$ is a decreasing sequence with limit greater or equal to $2\vare_0 >0$, so $\phi(\omega\setminus n)\geq 2\vare_0 >0$ for any $n\in\omega$. Therefore, for each $m,n\in\omega$, there exists $k>n$ such that $\phi(k\setminus n)>\vare_m$. Let $\left<k_i\right>_{i\in\omega}$, be an increasing sequence such that $k_{0}=0$ and $\phi(k_{i+1}\setminus \varp(x)(i))>\vare_{i}$, $i\in\omega$. Set $f_{j}(x)=\vare_{i}$ if $k_{i}\leq j<k_{i+1}$. Then $f_m(x)\ge\varepsilon_n$ if and only if $m<k_{n+1}$. Therefore, if $\phi(\{m\in\omega\colon f_{m}(x)\geq \vare_{n}\}\setminus k)<\vare_{n}$, then $k\geq\varp(x)(n)$, so $(o_\phi F)(x)(n)\geq \varp(x)(n)$ for any $n\in\omega$. 

This proves that $o$ is a cofinal function. Therefore by Lemma~\ref{eqlem}, the property $H^\Leftarrow(\oF_{\to_I},\twoheadrightarrow_I)$ holds. \hfill $\square$

\begin{cor}
Assume that $\non (\NN)=\cont$, and that there exists a $\cont$-Lusin set. Let $I$ be any analytic $P$-ideal. Then there exists $F=\left<f_n\right>_{n\in\omega}$, $f_{n}\colon [0,1]\to [0,1]$ for $n\in\omega$ with $f_{n}\to_{I} 0$ and $\vare>0$ such that for every $A\subseteq [0,1]$ with $m^{*}(A)\geq 1-\vare$,  $f_{n}\not\twoheadrightarrow_{I} 0$ on $A$ (the negation of the ideal version of the generalized Egorov's statement between equi-ideal and pointwise ideal convergence for analytic $P$-ideals is consistent with ZFC).
\end{cor}

Proof: We use Theorem \ref{thneg} and Lemma \ref{organicz}. \hfill $\square$

\section{Countably generated ideals}

Recall that an ideal $I$ over $\omega$  is countably generated (satisfies the chain condition) if there exists a sequence $\left<C_{i}\right>_{i\in\omega}$ of elements of $I$ such that $C_{i}\subseteq C_{i+1}$ for all $i\in\omega$ and for every $A\in I$, there exists $k\in\omega$ such that $A\subseteq C_{k}$. 

Let $\left<f_{n}\right>_{n\in\omega}$, $f_{n}\colon [0,1]\to[0,1]$, and let $I$ be an ideal on $\omega$. Recall the classic notion of ideal convergence on $A\subseteq [0,1]$:
\begin{description}
\item[pointwise ideal, $f_{n}\to_{I} 0$]
if and only if  $\forall_{\vare>0}\forall_{x\in A} \{n\in\omega\colon f_{n}(x)\geq\vare\}\in I$,
\item[quasinormal ideal, $f_n\xrightarrow{QN}_{I} 0$]
if and only if there exists a sequence of positive reals $\left<\vare_n\right>_{n\in\omega}$ such that $\vare_n\to_{I} 0$ and $\forall_{x\in A} \{n\in\omega\colon f_{n}(x)\geq \vare_n\}\in I$,
\item[uniform ideal, $f_{n}\rightrightarrows_{I} 0$] if and only if \[\forall_{\vare>0}\exists_{B\in I}\forall_{x\in A}\{n\in\omega\colon f_{n}(x)\geq\vare\}\subseteq B.\] 
\end{description}

The quasinormal convergence with respect to an ideal $I$ is also sometimes called $I$-equal convergence. Notice that in the case of countably generated ideals the generalized Egorov's statement holds between uniform ideal and quasinormal ideal convergence (see \cite[Theorem 3.2]{pdsdsk:iiecett}).

Let us therefore compare the pointwise and uniform ideal convergences. First, we show that the classic version (for measurable functions) of Egorov's Theorem holds in the case of convergence with respect to a countably generated ideal.

\begin{theo}
If $I\subseteq 2^{\omega}$ is a countably generated ideal and $f_{n}\colon [0,1]\to[0,1]$, $n\in\omega$ are measurable functions such that $f_{n}\to_I 0$ and $\vare>0$, then there exists a measurable set $B\subseteq [0,1]$ such that $m(B)\leq\vare$ and $f_{n}\rightrightarrows_I 0$ on $[0,1]\setminus B$.
\end{theo}

Proof: Assume that $I$ is countably generated and fix sets $\left<C_{i}\right>_{i\in\omega}$ such that $C_{i}\subseteq C_{i+1}$ for all $i\in\omega$ and for every $A\in I$, there exists $k\in\omega$ such that $A\subseteq C_{k}$. For $n,k\in\omega$, let
\[E_{n,k}=\left\{x\in[0,1]\colon \left\{m\in\omega\colon f_m(x)>\frac{1}{2^k}\right\}\setminus C_n\neq\zbp\right\}.\]
Notice that \[E_{n,k}=\bigcup_{m\in\omega\setminus C_n}\left\{x\in[0,1]\colon f_m(x)>\frac{1}{2^k}\right\}\] is measurable for each $n,k\in\omega$. Moreover, $E_{n+1,k}\subseteq E_{n,k}$ and $\bigcap_{n\in\omega} E_{n,k}=\zbp$ for all $k\in\omega$. Let $\vare>0$. For each $k\in\omega$, there exists $n_k\in\omega$ such that \[m(E_{n_k,k})\leq\frac{\vare}{2^{k+1}}.\] Let $B=\bigcup_{k\in\omega} E_{n_k,k}$. So $m(B)\leq\varepsilon$, and if $x\notin B$, then \[\left\{m\in\omega\colon f_m(x)>\frac{1}{2^k}\right\}\subseteq C_{n_k},\] for any $k\in\omega$, so $f_n\rightrightarrows_I 0$ on $[0,1]\setminus B$. \hfill $\square$

Let us consider the generalized Egorov's statement in this setting. The results presented below was proved by Joanna Jureczko using the method of T. Weiss (see \cite{tw:nget}) directly. We continue to apply the generalization of Pinciroli's method as presented above.

Assume that $I$ is countably generated, and fix sets $\left<C_{i}\right>_{i\in\omega}$ such that $C_{i}\subseteq C_{i+1}$ for all $i\in\omega$ and for every $A\in I$, there exists $k\in\omega$ such that $A\subseteq C_{k}$. We can assume that $C_{i+1}\setminus C_{i}\neq\zbp$ for all $i\in\omega$.

If $F=\left<f_{n}\right>_{n\in\omega}$, $f_{n}\to_{I} 0$, we define
\[(o_{\left<C_i\right>} F)(x) (n)=\min \left\{k\in\omega\colon \left\{m\in\omega\colon f_m(x)>\frac{1}{2^n}\right\}\subseteq C_k\right\}.\]

Notice that if $A\subseteq [0,1]$, then $f_n\rightrightarrows_I 0$ on $A$ if and only if $(o_{\left<C_i\right>} F)[A]$ is bounded, and so $H^\Rightarrow(\oF_{\to_I},\rightrightarrows_I)$ holds.  Therefore, we get the following theorem.

\begin{cor}
Assume that $\non(\NN)<\bont$. Let $I$ be any countably generated ideal, and let $\vare>0$. Let $F=\left<f_{n}\right>_{n\in\omega}$, $f_{n}\colon [0,1]\to [0,1]$, for $n\in\omega$  be such that $f_{n}\to_{I} 0$. Then there exists $A\subseteq [0,1]$ with $m^{*}(A)\geq 1-\vare$ such that $f_{n}\rightrightarrows_{I} 0$ on $A$ (the ideal version of the generalized Egorov's statement between uniform ideal and pointwise ideal convergence for countably generated ideals is consistent with ZFC).
\end{cor}

Proof: Apply Theorem \ref{thpos}. \hfill $\square$

\begin{lemma}\label{ogranccc}
For any $\varp\colon [0,1]\to \omega^{\omega}$ there exists $F=\left<f_{n}\right>_{n\in\omega}$, $f_{n}\colon [0,1]\to [0,1]$, $f_{n}\to_{I} 0$ for $n\in\omega$ such that $o_{\left<C_i\right>} F= \varp$. In particular, $H^\Leftarrow(\oF_{\to_I},\rightrightarrows_I)$ holds.
\end{lemma}
Proof: Without a loss of generality we can assume that $\varp(x)$ is increasing for all $x\in [0,1]$. Let $x\in [0,1]$. Let $f_{j}(x)=1/2^n$ if and only if $j\in C_{\varp(x)(n+1)}\setminus C_{\varp(x)(n)}$. \hfill $\square$

\begin{cor}
Assume that $\non (\NN)=\cont$, and that there exists a $\cont$-Lusin set. Let $I$ be any countably generated ideal. Then there exists $F=\left<f_n\right>_{n\in\omega}$, $f_{n}\colon [0,1]\to [0,1]$ for $n\in\omega$ with $f_{n}\to_{I} 0$, and $\vare>0$ such that for all $A\subseteq [0,1]$ with $m^{*}(A)\geq 1-\vare$,  $f_{n}\not\rightrightarrows_{I} 0$ on $A$ (the negation of the ideal version of the generalized Egorov's statement between uniform ideal and pointwise ideals convergence for countably generated ideal is consistent with ZFC).
\end{cor}

Proof: Apply Theorem \ref{thneg} and Lemma \ref{ogranccc}. \hfill $\square$

\section{$I^*$ convergence for countably generated ideals}

As before, let $\left<f_{n}\right>_{n\in\omega}$, $f_{n}\colon [0,1]\to[0,1]$, and let $I$ be an ideal on $\omega$. We have the following notion of convergence $A\subseteq [0,1]$ (see \cite{pdsdsk:iiecett}):
\begin{description}
\item[$I^*$-pointwise, $f_{n}\to_{I^*} 0$] if and only if for all $x\in A$, there exists $M=\{m_{i}\colon i\in\omega\}\subseteq\omega$, $m_{i+1}>m_i$ for $i\in\omega$ such that $\omega\setminus M\in I$ and $f_{m_i}(x)\to 0$,
\item[$I^*$-quasinormal, $f_n\xrightarrow{QN}_{I^*} 0$]
if and only if there exists $M=\{m_{i}\colon i\in\omega\}\subseteq\omega$, $m_{i+1}>m_i$ for $i\in\omega$ such that $\omega\setminus M\in I$ and $f_{m_i}\xrightarrow{QN} 0$ on $A$,
\item[$I^*$-uniform, $f_{n}\rightrightarrows_{I^*} 0$] if and only if there exists $M=\{m_{i}\colon i\in\omega\}\subseteq\omega$, $m_{i+1}>m_i$ for $i\in\omega$ such that $\omega\setminus M\in I$ and $f_{m_i}\rightrightarrows 0$ on $A$. 
\end{description}

Notice that for any ideal $I$, the generalized Egorov's statement holds between $I^*$-uniform  and $I^*$-quasinormal convergence (see \cite[Theorem 3.3]{pdsdsk:iiecett}).

Let us therefore compare the pointwise and uniform ideal convergences.
First, we show that the classic version (for measurable functions) of Egorov's Theorem holds in the case of $I^*$-convergence with respect to a countably generated ideal $I$.

\begin{theo}
If $I\subseteq 2^{\omega}$ is a countably generated ideal and $f_{n}\colon [0,1]\to[0,1]$, $n\in\omega$ are measurable functions such that $f_{n}\to_{I^*} 0$ and $\vare>0$, then there exists a measurable set $B\subseteq [0,1]$ such that $m(B)\leq\vare$ and $f_{n}\rightrightarrows_{I^*} 0$ on $[0,1]\setminus B$.
\end{theo}

Proof: Assume that $I$ is countably generated and fix $\left<C_n\right>_{n\in\omega}$ such that for all $A\in I$, there exists $n\in\omega$ with $A\subseteq C_n$. Let $\omega\setminus C_{n}=\{m_{i,n}\colon i\in\omega\}$, $m_{i+1,n}>m_{i,n}$, $i,n\in\omega$, and 
\[F_n=\left\{x\in[0,1]\colon \lim_{i\in\omega}f_{m_{i,n}}(x)=0\right\}\]
Obviously, $F_{n}\subseteq F_{n+1}$ for $n\in\omega$ and $\bigcup_{n\in\omega} F_n=[0,1]$. Moreover, \[F_n=\bigcap_{i\in\omega}\bigcup_{j\in\omega}\bigcap_{k\geq j}\left\{x\in [0,1]\colon f_{m_{k,n}}(x)<\frac{1}{2^i}\right\}\] is measurable. Therefore, there exists $N\in\omega$ such that $m(F_N)\geq 1-\vare/2$. Now apply the classic Egorov's Theorem for the set $F_N$, $\left<f_{m_{i,N}}\right>_{i\in\omega}$ and $\vare/2$ to get a set $A\subseteq F_N$ such that $f_{m_{i,N}}$ converges uniformly on $F_N\setminus A$ and $m(A)<\vare/2$. Let $B=A\cup([0,1]\setminus F_N)$. We get that $f_n\rightrightarrows_{I^*} 0$ on $[0,1]\setminus B$ and $m(B)\leq\vare$. \hfill $\square$  

Let us consider the generalized Egorov's statement in this setting. Assume that $I$ is countably generated  and fix $\left<C_n\right>_{n\in\omega}$ such that for all $A\in I$, there exists $n\in\omega$ such that $A\subseteq C_n$. Let $F=\left<f_n\right>_{n\in\omega}$ be such that $f_{n}\to_{I^*} 0$. 
Let $F=\langle f_n\rangle_{n\in\omega}$ be such that $f_n\to_{I^*}0$.
For $x\in[0,1]$ define $o_{\langle C_i\rangle}(F)(x)=\psi\in\omega^\omega$ by
\begin{align*}
&\psi(0)=\min\left\{n\in\omega:
\langle f_m\rangle_{m\in\omega\setminus C_n}\to 0\right\},\\
&\psi(n)=\min\left\{m\in\omega:
\forall_{\substack{l\in\omega\setminus C_{\psi(0)}\\l>m}}\ f_l(x)<\frac{1}{2^{n}}\right\},
\quad n>0.
\end{align*}

Obviously, $o_{\left<C_i\right>} F$ is bounded if and only if $f_n\rightrightarrows_{I^*} 0$, and so the property $H^\Rightarrow(\oF_{\to_{I^*}},\rightrightarrows_{I^*})$ holds.

Therefore, we get the following theorem.

\begin{cor}
Assume that $\non(\NN)<\bont$. Let $I$ be any countably generated ideal, and let $\vare>0$ and $F=\left<f_{n}\right>_{n\in\omega}$, $f_{n}\colon [0,1]\to [0,1]$ for $n\in\omega$, with $f_{n}\to_{I^*} 0$. Then there exists $A\subseteq [0,1]$ with $m^{*}(A)\geq 1-\vare$ such that $f_{n}\rightrightarrows_{I^*} 0$ on $A$ (the ideal version of the generalized Egorov's statement between uniform $I^*$ and pointwise $I^*$ convergence for countably generated ideals is consistent with ZFC).
\end{cor}

Proof: Apply Theorem~\ref{thpos}. \hfill $\square$

\begin{lemma}\label{ograngw}
For any $\varp\colon [0,1]\to \omega^{\omega}$, there exist $F=\left<f_{n}\right>_{n\in\omega}$, $f_{n}\colon [0,1]\to [0,1]$, $f_{n}\to_{I^*} 0$ for $n\in\omega$ such that $o_{\left<C_i\right>} F\geq \varp$. In particular, the condition $H^\Leftarrow(\oF_{\to_{I^*}},\rightrightarrows_{I^*})$ holds.
\end{lemma}
Proof: It is enough to prove the lemma for $\varphi$ such that $\varphi(x)$ is
increasing for all $x\in[0,1]$. Let $x\in [0,1]$. Let $\omega\setminus C_{\varp(x)(0)}=\{m_{i}\colon i\in\omega\}$, $m_{i+1}>m_i$ for $i\in\omega$. Let $f_{j}(x)=1$ for $j\in C_{\varp(x)(0)}$ and let $f_j(x)=1/2^n$ if $j\in (\omega\setminus C_{\varp(x)(0)})\cap\{i\in\omega\colon \varp(x)(n)\leq i <\varp(x)(n+1)\}$. \hfill $\square$

\begin{cor}
Assume that $\non (\NN)=\cont$, and that there exists a $\cont$-Lusin set. Let $I$ be any countably generated  ideal. Then there exists $F=\left<f_n\right>_{n\in\omega}$, $f_{n}\colon [0,1]\to [0,1]$ for $n\in\omega$, with $f_{n}\to_{I^*} 0$, and $\vare>0$ such that for all $A\subseteq [0,1]$ with $m^{*}(A)\geq 1-\vare$,  $f_{n}\not\rightrightarrows_{I^*} 0$ on $A$ (the negation of the ideal version of the generalized Egorov's statement between uniform $I^*$ and pointwise $I^*$ convergence for countably generated ideals is consistent with ZFC).
\end{cor}

Proof: Apply Theorem~\ref{thneg} and Lemma \ref{ograngw}. \hfill $\square$

\section{Ideals $\Fin^{\alpha}$}

Given an ideal $I\subseteq \omega$ and a sequence $\left<I_n\right>_{n\in\omega}$ of ideals of $\omega$, we can consider an ideal $I\mhyphen\prod_{n\in\omega} I_n$ on $\omega^2$ called the $I$-product of the sequence of ideals $\left<I_n\right>_{n\in\omega}$ and define it in the following way. For any $A\subseteq \omega^2$,
\[A\in I\mhyphen\prod_{n\in\omega} I_n \leftrightarrow \{n\in\omega\colon A_{(n)}\notin I_n\}\in I,\]
where $A_{(n)}=\{m\in\omega\colon \left<n,m\right>\in A\}$ (see \cite{nm:zicf}). If $I_n=J$ for any $n\in\omega$, we usually denote $I\mhyphen\prod_{n\in\omega} I_n$ as $I\times J$.

Fix a bijection $b\colon \omega^2\to\omega$ and a bijection $a_{\beta}\colon \omega\setminus \{0\}\to \beta$ for any limit $\beta<\omega_1$. The ideals $\Fin^\alpha$, $\alpha<\omega_1$, are defined inductively (see \cite{nm:zicf}) in the following way. Let $\Fin^1=\Fin$ be the ideal of finite subsets of $\omega$. We set $\Fin^{\alpha+1}=\{b[A]\colon A\in \Fin\times \Fin^\alpha\}$ and for limit $\beta<\omega_1$, let $\Fin^{\beta}=\{b[A]\colon A\in \Fin\mhyphen\prod_{i\in\omega} \Fin^{a_{\beta}(i+1)}\}$. 

In \cite[Theorem 3.25]{nm:zicf}, N.~Mro\.{z}ek proves that ideal $\Fin^{\alpha}$ for any $\alpha<\omega_1$ satisfies the Egorov's theorem for ideals (between uniform ideal and pointwise ideal convergences). 

Let $\oF_\alpha=\oF_{\to_{\Fin^{\alpha}}}$. We get the following theorem.

\begin{theo}
Assume that $\non(\NN)<\bont$. Let $0<\alpha<\omega_1$, and let $\vare>0$ and $F=\left<f_{n}\right>_{n\in\omega}$, $f_{n}\colon [0,1]\to [0,1]$ for $n\in\omega$, with $f_{n}\to_{\Fin^{\alpha}} 0$. Then there exists $A\subseteq [0,1]$ with $m^{*}(A)\geq 1-\vare$ such that $f_{n}\rightrightarrows_{\Fin^{\alpha}} 0$ on $A$ (the ideal version of the generalized Egorov's statement between uniform $\Fin^{\alpha}$ and pointwise $\Fin^{\alpha}$ convergence is consistent with ZFC).
\end{theo}

Proof: We define $o_\alpha\colon \oF_{\alpha}\to (\omega^{\omega})^{[0,1]}$ in the following way. Let $\varepsilon_n=\frac{1}{2^n}$ for $n\in\omega$, and let 
\[\oF_\alpha^n=\{\langle f_n\rangle_{n\in\omega}:\forall_{n\in\omega} f_n:[0,1]\to[0,1]\land \forall_{x\in [0,1]}\{q\in\omega\colon f_q(x)\geq \vare_n\}\in\Fin^\alpha\}.\] 
First, define $o^n_\alpha \colon \oF_\alpha^n\to (\omega^{\omega})^{[0,1]}$, $n\in \omega, 0<\alpha<\omega_1$, by induction on $\alpha$. Let 
\[M_{1, n, x}=\min\{p\in\omega\colon\forall_{q\geq p} f_q(x)<\varepsilon_n\},\]
and let
\[(o^n_1 F)(x)(k)=M_{1, n, x}\] 
be a constant sequence. Given $o^n_\alpha$, let
\[M_{\alpha+1, n, x}=\min\left\{p\in\omega\colon \forall_{q\geq p} \{m\in\omega\colon f_{b(q,m)}(x)\geq \varepsilon_n\}\in \Fin^{\alpha}\right\},\]
and
\[
(o^n_{\alpha+1} F)(x)(k)=\begin{cases} 
M_{\alpha+1, n, x}& \text{ for } k=b(p,q),\\ &  p<M_{\alpha+1, n, x}+1, q\in\omega,\\
(o^n_{\alpha}\left<f_{b(p-1,r)}\right>_{r\in\omega}) (x)(q)& \text{ for } k=b(p,q),\\  & p\geq M_{\alpha+1, n, x}+1, q\in\omega.\end{cases}
\]
This definition is correct, since $\left<f_{b(p-1,r)}\right>_{r\in\omega}\in\oF_\alpha^n$ for $p\geq M_{\alpha+1, n, x}+1$. 

Moreover, for limit $\beta<\omega_1$, let
\[M_{\beta, n, x}=\min\left\{p\in\omega\colon \forall_{q\geq p} \{m\in\omega\colon f_{b(q,m)}(x)\geq \varepsilon_n\}\in \Fin_{a_\beta(q)}\right\}\]
and
\[(o^n_{\beta} F)(x)(k)=\begin{cases} 
M_{\beta, n, x}& \text{ for } k=b(p,q),\\ & p<M_{\beta, n, x}+1, q\in\omega,\\
(o^n_{a_{\beta}(p-1)}\left<f_{b(p-1,r)}\right>_{r\in\omega})(x)(q)& \text{ for } k=b(p,q),\\ & p\geq M_{\beta, n, x}+1, q\in\omega.\end{cases}
\]
This definition is correct, since, for each $p\geq M_{\beta, n, x}+1$, $\left<f_{b(p-1,r)}\right>_{r\in\omega}\in\oF_{a_{\beta}(p-1)}^n$.

Notice that $\oF_{\alpha}\subseteq \oF_{\alpha}^n$, for any $n\in\omega$. Therefore, finally let
\[(o_\alpha F)(x)(k)= (o^n_\alpha F)(x)(m),\]
for $k=b(n,m)$, $n,m\in\omega$.

Now, notice that if $F=\left<f_r\right>_{r\in\omega}\in \oF_\alpha$, and $o_\alpha F$ is bounded on a set $A\subseteq [0,1]$, then $f_r\rightrightarrows_{\Fin^\alpha} 0$ on $A$. Indeed, if $o_\alpha F$ is bounded, then for each $n\in\omega$, $o^n_\alpha F$ is bounded. If so, $\{m\in\omega\colon \sup_{x\in A} f_m(x)\geq \varepsilon_n\}\in \Fin^\alpha$, for all $n\in\omega$. We fix $n\in\omega$ and prove this statement by induction on $\alpha<\omega_1$. Let $(o^n_\alpha F)(x)(k)< a_{k,n}$ for all $x\in A, k\in\omega$ and some $\left<a_{k,n}\right>_{k\in\omega}\in \omega^\omega$. If $\alpha=1$, we get $f_q(x)<\varepsilon_n$ for all $x\in A$ and all $q\geq a_0$, so $\{m\in\omega\colon \sup_{x\in A} f_m(x)\geq \varepsilon_n\}\in \Fin$. Now, assume that the statement holds for some $\alpha<\omega_1$. Then for all $x\in A$, $M_{\alpha+1, n, x}<a_{b(0,0)}$, so for all $p\geq a_{b(0,0)}$, $o^n_{\alpha}\left<f_{b(p-1,r)}\right>_{r\in\omega} $ is bounded by $\left<a_{b(p,q)}\right>_{q\in\omega}$, and thus by the induction hypothesis, $\{r\in\omega\colon \sup_{x\in A} f_{b(p-1,r)}\geq\varepsilon_n\}\in\Fin^\alpha$ for all $p\geq a_{b(0,0)}$. Therefore, $\{m\in\omega\colon \sup_{x\in A} f_m(x)\geq \varepsilon_n\}\in \Fin^{\alpha+1}$. Analogous reasoning can be easily applied for limit $\beta<\omega_1$. This proves that $H^\Rightarrow(\oF_{\alpha},\rightrightarrows_{\Fin^{\alpha}})$ holds.

Therefore, by Theorem \ref{thpos}, there exists $A\subseteq [0,1]$ with $m^{*}(A)\geq 1-\vare$ such that $f_{n}\rightrightarrows_{\Fin^{\alpha}} 0$ on $A$. \hfill $\square$

\begin{theo}
Assume that $\non (\NN)=\cont$, and that there exists a $\cont$-Lusin set. Let $0< \alpha<\omega_1$. Then there exist $\left<f_{n}\right>_{n\in\omega}\in\oF_\alpha$ and $\vare>0$ such that for all $A\subseteq [0,1]$ with $m^{*}(A)\geq 1-\vare$, $f_{n}\not\rightrightarrows_{\Fin^\alpha} 0$ on $A$ (the negation of the ideal version of the generalized Egorov's statement between uniform $\Fin^{\alpha}$ and pointwise $\Fin^{\alpha}$ convergence for countably generated ideals is consistent with ZFC).
\end{theo}

Proof: As before, let $\varepsilon_n=1/2^n$, $n\in\omega$. This time, we define $o_\alpha$ in a different way then in the previous proof. Namely, let
\[(o_{\alpha} F)(x)(n)= M_{\alpha,n,x},\]
where $M_{\alpha,n,x}$ is defined as in the previous proof. Notice that if $F=\left<f_n\right>_{n\in\omega}$ is such that $f_n\rightrightarrows_{\Fin^\alpha} 0$ on a set $A\subseteq [0,1]$, then $\{m\in\omega\colon \sup_{x\in A} f_m(x)\geq \varepsilon_n\}\in \Fin^\alpha$ for all $n\in\omega$. If $\alpha=1$, this means that $\min\{p\in\omega\colon\forall_{q\geq p} f_q(x)<\varepsilon_n\}=M_{1,n,x}=o_{1} F(x)(n)$ is bounded on $A$. If $\alpha$ is a limit ordinal, then for all $n\in\omega$, there exists $M_n$ such that for all $q\geq M_n$,  $\{m\in\omega\colon f_{b(q,m)}(x)\geq \varepsilon_n\}\in \Fin_{a_\alpha(q)}$. In other words, $M_{\alpha,n,x}=o_{\alpha} F(x)$ is bounded on $A$. Similar argument can be used in the case of a successor ordinal $\alpha>1$.

Moreover, fix any $\varp\colon [0,1]\to \omega^{\omega}$. Without a loss of generality, assume that for $x\in [0,1]$, $\varp(x)$ is increasing. There exists $F=\left<f_n\right>_{n\in\omega}\in\oF$ such that $o_{\alpha}(F)\geq\varp$. It is obvious for $\alpha=1$. For $\alpha>1$, let $f_n(x)=\varepsilon_k$ for $k= b(i,j)$, $\varp(x)(k)\leq n< \varp(x)(k+1)$. Therefore $H^\Leftarrow(\oF_{\alpha},\rightrightarrows_{\Fin^{\alpha}})$ holds.

In conclusion, by Theorem \ref{thneg}, there exist $\left<f_{n}\right>_{n\in\omega}\in\oF$ and $\vare>0$ such that for all $A\subseteq [0,1]$ with $m^{*}(A)\geq 1-\vare$, $f_{n}\not\rightrightarrows_{\Fin^\alpha} 0$ on $A$.  \hfill $\square$

\medskip
\noindent {\bf Acknowledgment}. The author is very grateful to the referee for a number of helpful suggestions for improvement of the paper.


\ifx\undefined\bysame
\newcommand{\bysame}{\leavevmode\hbox to3em{\hrulefill}\,}
\fi

\end{document}